\documentclass[12pt]{amsart}
\usepackage{amsmath,amssymb,amscd,amsfonts}
\newtheorem{prop}{Proposition}[section]
\newtheorem{thm}[prop]{Theorem}

\def \R {{\mathbb R}}

        \def\zed{{\mathbb Z}}
        \def\Z{{\zed}}
        
\input epsf.tex
\begin{document}

\title{Tilings, Tiling Spaces and Topology}

\author{Lorenzo Sadun}

\address{Lorenzo Sadun: Department of Mathematics, The University of
  Texas at Austin, Austin, TX 78712-1082 U.S.A.}
\email{sadun@math.utexas.edu}



\begin{abstract}
  To understand an aperiodic tiling (or a quasicrystal modeled on an
  aperiodic tiling), we construct a space of similar tilings, on which
  the group of translations acts naturally.
  This space is then an (abstract) dynamical system. Dynamical
  properties of the space (such as mixing, or the spectrum of the
  translation operator) are closely related to bulk properties of
  individual tilings (such as the diffraction pattern). The topology
  of the space of tilings, particularly the Cech cohomology, gives
  information on how original tiling may be deformed. Tiling spaces
  can be constructed as inverse limits of branched manifolds. 
\end{abstract}

\maketitle

\markboth{Sadun}{Tilings, Tiling Spaces and Topology}
\newpage

\section{Physical and mathematical questions}

Tilings are studied by mathematicians, physicists and material scientists.
To a physicist, an aperiodic tiling, such as the Penrose tiling,
is a model for a form of matter that is neither crystaline nor disordered. 
The most interesting questions are then about the physical properties of
the material being modeled:
\begin{itemize}

\begin{item}{P1.} What is the x-ray diffraction pattern of the material?
This is equivalent to the Fourier transform of the autocorrelation function
of the positions of the atoms.  Sharp peaks are the hallmark of ordered 
materials, such as crystals and quasicrystals.  
\end{item}

\begin{item} P2. What are the possible energy levels of electrons in the
material?  The locations of the atoms determine a quasiperiodic potential, 
and the spectrum of the corresponding Schr\"odinger Hamiltonian has infinitely
many gaps.  What are the energies of these gaps, and what is the density 
of states corresponding to each gap?
\end{item}

\begin{item} P3. Can you really tell the internal structure of the material
from diffraction data?  What deformations (either local or non-local)
of the molecular structure are consistent with the combinatorics of the 
molecular bonds?  Which of these are detectable from diffraction data?
\end{item}
\end{itemize}

\begin{figure}
\vbox{\epsfxsize=2truein\epsfbox{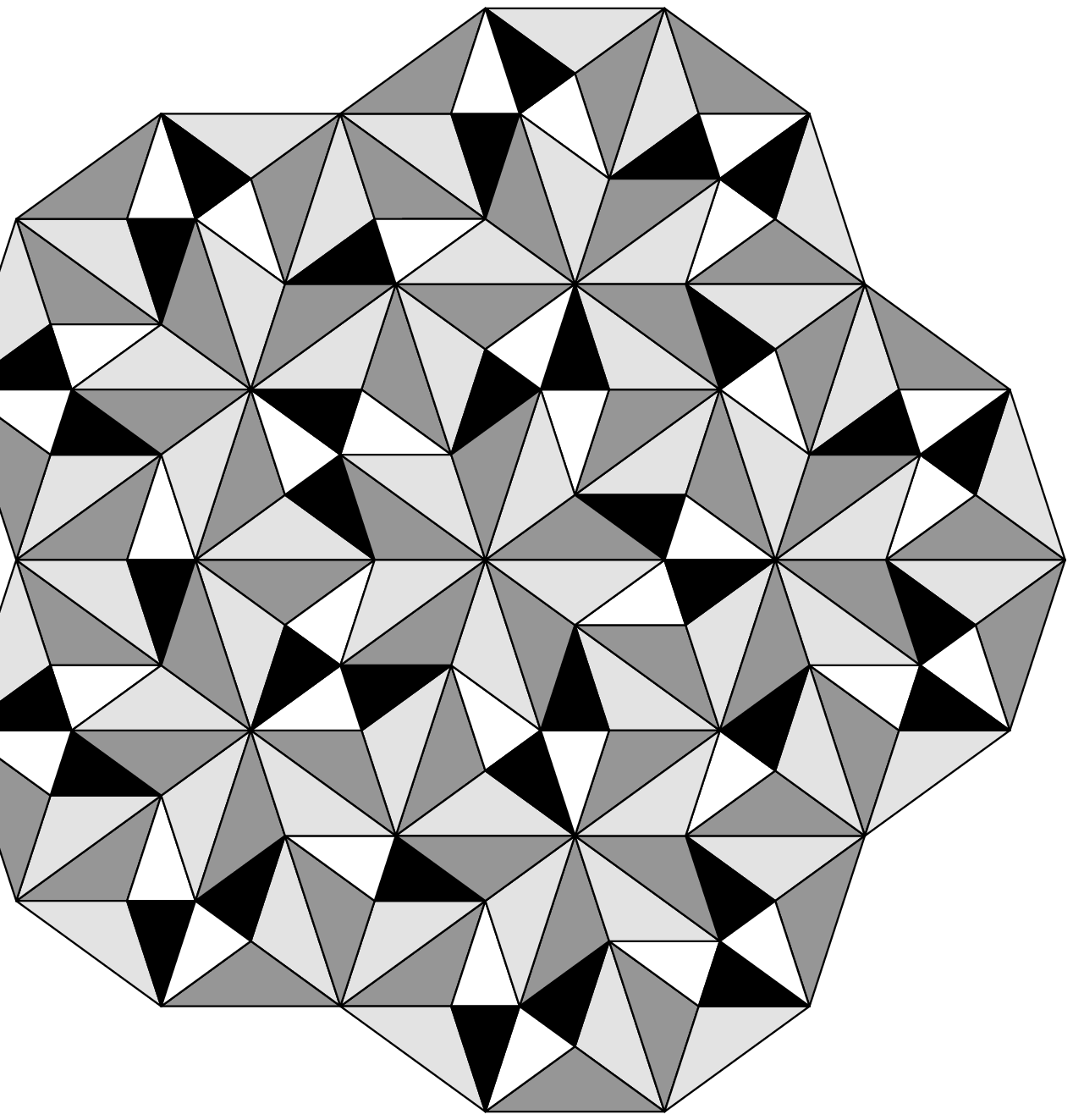}
\hskip 0.8 in \epsfxsize=2truein\epsfbox{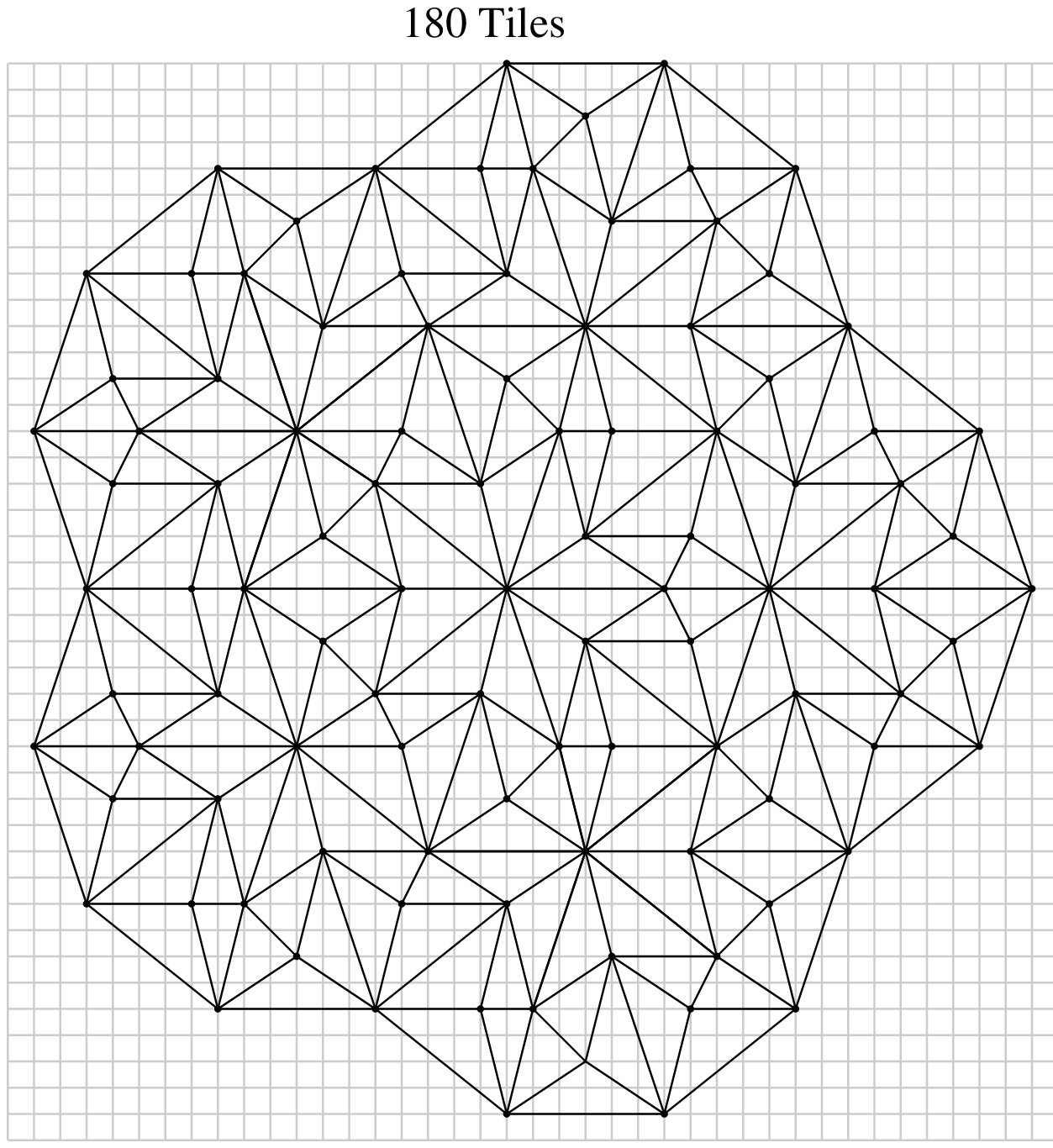}}
\caption{Patches of a Penrose tiling and a deformed Penrose tiling}
\end{figure}


As an example of P3, consider the tilings in figure 1.  The first is a
piece of the standard Penrose (rhomb) tiling, with each rhomb cut
into two triangles. The second tiling is combinatorially identical to
the first, but the shapes and sizes of the (40) different species of
tiles have been changed.  In the first tiling, the diffraction pattern
has rotational symmetry, since the underlying tiling has symmetry.  In 
the second tiling, the vertices all lie on lattice points, so the diffraction
pattern is periodic with period $(2 \pi \Z)^2$. However, we will see that
the diffraction peaks of the two tilings are related by a simple linear
transformation.  In other words, the first tiling and a linear transformation
of the second have very different geometry but qualitatively similar
diffraction patterns.

A mathematician takes a different approach to tilings.  To a
mathematician, a tiling $T$ is a point in a tiling space $X_T$.  The
other points of $X_T$ are tilings with the same properties as $T$ (see
the construction below).  If $T$ has nice properties (finite local
complexity, repetitivity, well-defined patch frequencies), then $X_T$
will have nice properties (compactness, minimality as a dynamical
system, unique ergodicity), and we can ask the following mathematical
questions:

\begin{itemize}
\begin{item} M1. What is the topology of $X_T$?  What does the neighbourhood
of a point of $X_T$ look like?  What are the (\v Cech) cohomology groups of 
$X_T$?
\end{item}

\begin{item} M2. There is a natural action of the group $\R^d$ 
of translations on $X_T$. 
This makes $X_T$ into a dynamical system, with $d$ commuting flows.  
What are the ergodic measures on $X_T$?  For each such measure, what is
the spectrum of the generator of translations (think: momentum operator) on
$L^2(X_T)$?  This is called the {\em dynamical spectrum} of $X_T$. 
\end{item}

\begin{item} M3. From the action of the translation group on $X_T$, one can
construct a $C^*$ algebra.  What is the K-theory of this $C^*$-algebra?
\end{item}
\end{itemize}

Remarkably, each math question about $X_T$ answers a physics question about
a material modeled on $T$. M1 answers P3, M2 answers P1, and M3 (in large part)
answers P2.  In this paper we will review the construction of $X_T$ and
explain the connection between M1 and P3,
relying on \cite{SW, CS, inverse} for details and proofs.  The connection
between dynamical and diffraction spectrum is an old story, and we refer 
the reader to \cite{Dw, Hof, BL, Goure, LMS}.  For the relation
between K-theory and gap-labeling (i.e., the density of states associated
to gaps) see \cite{Bel, BHZ, BBG}. 
 
\section{The space $X_T$}

A point in our tiling space will be a tiling of $\R^d$ (for
simplicity, we'll stick with $d=2$, but the construction is the same
in all dimensions).  Note that $\R^2$ isn't just a plane --- it is a
plane with a distinguished point, namely the origin.  Translating the
tiling to the right is equivalent to moving the origin to the left,
{\em and results in a different tiling}.  The orbit of $T$ is the set of
all translates of $T$.  

We put a metric on tilings as follows.  We say two tilings $T_1$ and
$T_2$ are $\epsilon$-close if they agree on a ball of radius
$1/\epsilon$ around the origin, up to a further translation by
$\epsilon$ or less. If $T_1=T-x$ and $T_2=T-y$ are translates of $T$,
this implies that a ball of radius $1/\epsilon$ around $x$ in $T$
looks just like a ball of radius $1/\epsilon$ around $y$.  The space 
$X_T$ is the completion of the orbit of $T$ in this metric.  Put another
way, it is the set of tilings $S$ with the property that every patch of
$S$ can be found somewhere in $T$.  Other names for $X_T$ are the
{\em continuous hull of $T$} and the {\em local isomorphism class of $T$}. 

The first topological question is whether $X_T$ is compact.  If $T$ has
an infinite number of tile types, or if the tiles can fit together in an
infinite variety of ways (e.g., a continuous shear along a common edge), then
it is easy to construct a sequence of tilings whose behavior at the origin
never settles down, and which therefore does not have a convergent subsequence.
However, if $T$ has only a {\em finite} number of tile types, and these tiles 
only fit together in only a finite number of ways, then there are only a finite
number of patterns (up to translation) that can appear in a fixed region.  
This is called {\em finite local complexity}, or FLC.  If $T$ is an FLC tiling,
then any sequence of translates of $T$ has a subsequence that converges on 
a ball of radius 1 around the origin.  That sequence has a subsequence that
converges on a ball of radius 2, etc.  Using a Cantor diagonalization 
argument, we can easily produce a subsequence that converges everywhere. In 
other words, we have proved that

\begin{thm}
The tiling space $X_T$ is compact if and only if $T$ has finite local
complexity. 
\end{thm}

Every FLC tiling space is MLD-equivalent (see below)
to a tiling space where the tiles are polygons 
that meet full-edge to full-edge.  (See \cite{Pr} for a discussion of
the derived Voronoi tilings that accomplish this equivalence.)
Without loss of generality, then, 
we can restrict our attention to tilings of this sort.

A tiling $T$ is called
{\em repetitive} if, for every patch $P$ in $T$, there is a radius
$R_P$ such that every ball of radius $R_P$ contains at least one copy
of $P$.  This guarantees that every tiling $S \in X_T$ contains the
patch $P$, hence that the set of patches of $S$ is exactly the same as
the set of patches of $T$, and hence that $X_S = X_T$.  In other
words, there is nothing special about $T$ itself, and the entire
tiling space can be recovered from any tiling in the space.
Equivalently, the orbit of every tiling $S \in X_T$ is dense in $X_T$.
In dynamical systems language, the tiling dynamical system $X_T$ is
{\em minimal}.

If $T$ is a periodic tiling (e.g., a checkerboard), invariant under
translation by a lattice $L$, then it is easy to see that $X_T$ is the
torus $\R^2/L$.  When $X_T$ is totally non-periodic, meaning there are
{\em no} vectors $x \in \R^2$ for which $T = T+x$, then $X_T$ is much 
more complicated.  Also much more interesting.

Locally, $X_T$ looks like a disk crossed with a totally disconnected
set.  By definition, if $S \in X_T$, 
then an $\epsilon$-neighbourhood of $S$ is all
tilings that agree with a small translate of $S$ on a ball of radius
$1/\epsilon$.  This gives $2$ continuous degrees of freedom (the small
translations), and a number of discrete degrees of freedom (choices on
how to extend the tiling beyond the $1/\epsilon$ ball). If the
original tiling $T$ is repetitive and non-periodic, then there are
infinitely many discrete degrees of freedom, and the discrete choices
yield a Cantor set.

To somebody who is used to smooth manifolds, a set that is locally a
disk crossed with a Cantor set is bizarre. For instance $X_T$ is
connected (since it has a path-connected dense set, namely the orbit
of $T$) but not path-connected (since the path component of $T$ is
merely the orbit of $T$). However, such structures are common in
dynamical systems.

There are several important notions of equivalence for tiling spaces.
Let $S$ and $T$ be different tilings.  We say that $X_S$ and $X_T$ are
{\em topologically conjugate} if there is a continuous map 
$\rho: X_S \to X_T$, 
with a continuous inverse, that commutes with translations. As dynamical
systems, $X_T$ and $X_S$ then have identical properties.  In particular,
they have the same mixing properties and the same dynamical spectrum.

An even stronger equivalence is {\em mutual local derivability}, or MLD. 
We say that $X_S$ and $X_T$ are MLD if they
are topologically conjugate and the map $\rho$ is local.  That is, if there
exists a radius $R$ such that, when $T_1, T_2 \in X_T$ agree exactly on a
ball of radius $R$ around a point $x$, then $\rho(T_1)$ and $\rho(T_2)$
agree exactly on a ball of radius 1 around $x$.  We shall see that the two
tiling spaces of Figure 1 are topologically conjugate (up to a fixed linear
transformation), but not MLD.

\section{Inverse limits}

One way to get a handle on tiling spaces is via inverse limits.  If 
$K_1, K_2, \ldots$ are topological spaces and $\sigma_n: K_n \to K_{n-1}$
are continuous maps, then the {\em inverse limit space} 
$\lim_{\leftarrow} K_n$ is defined as 
$$ \lim_{\leftarrow} K_n = \{ (x_1,x_2,\ldots) \in \prod_n K_n| \sigma_n(x_n)=
x_{n-1}\} \eqno(1) $$
In other words, a point in $\lim_{\leftarrow} K_n$ is a point $x_1 \in K_1$,
together with a point $x_2 \in K_2$ that maps to $x_1$, together with a
point $x_3 \in K_3$ that maps to $x_2$, and so on. Two sequences are considered
close if their first $N$ terms are close, for $N$ large.  The space $K_n$ is
called the {\em $n$-th approximant} to $K_n$, $x_n$ determines
$x_{n-1}, x_{n-2},\ldots, x_1$, 

A simple example of an inverse limit 
is when each $K_n$ is the circle $\R/\Z$ and each
$\sigma_n$ is multiplication by 2.  This is called the {\em dyadic
  solenoid}.  For each $x_1 \in \R/\Z$, there are two possibilities
for $x_2$, for each of these there are two possibilities for $x_3$,
for each of these there are two possibilities for $x_4$, and so on.
The points within $\epsilon$ of $(x_1,x_2,\ldots)$ can differ by a
continuous motion (adding $\delta$ to $x_1$, $\delta/2$ to $x_2$,
$\delta/4$ to $x_3$, and so on), or by keeping the first $N$ terms the
same and then making arbitrary choices for $x_{N+1}$, $x_{N+2}$, etc.
This is exactly the same local structure as a 1-dimensional tiling
space, namely one continuous degree of freedom and infinitely many 
discrete possibilities.

A nice property of inverse limit spaces is that their (\v Cech) cohomologies
are easy to compute.  The cohomology of the {\em inverse} limit is the
{\em direct} limit of the cohomologies of the approximants: 
$$ H^k(\lim_{\leftarrow} K_n) = \lim_{\rightarrow} H^k(K_n) \eqno(2) $$
For the dyadic solenoid,
each $K_n$ has $H^0 = H^1 = \Z$.  The pullback map $\sigma_n^*$ is 
multiplication by 1 on $H^0$ and multiplication by 2 on $H^1$, so the dyadic
solenoid has $H^0=\Z$, $H^1=\Z[1/2]$.  That is, $H^1$ is 
the set of rational numbers whose denominators are powers of 2.  

Most computations of the cohomology of tiling spaces use inverse
limits and equation (2) \cite{AP, ORS, BBG, G}.  However, for certain
``cut-and-project'' tilings, other techniques have been developed \cite{FHK}.

\section{Tiling spaces are inverse limits}

Spaces of tilings can always be viewed as inverse limits. Anderson and
Putnam \cite{AP} first observed this for substitution tilings. 
Others extended the idea to larger classes of tiling space
\cite{ORS, BBG, BG}.  The following general construction is due to G\"ahler 
\cite{G}, and was further generalized in \cite{inverse}.

\begin{thm}\label{4.1} If $T$ is a tiling with finite local complexity, then
$X_T$ is the inverse limit of a sequence of compact branched
surfaces $K_1, K_2, \ldots$ and continuous maps $\sigma_n: K_n \to K_{n-1}$.
\end{thm}

A point in $K_n$ will be a set of instructions for placing a tile
containing the origin, a ring of tiles around it (the {\em first
  corona}), a second ring around that, and on out to the $n$-th
corona. The map $\sigma_n:K_n \to K_{n-1}$ simply forgets the
outermost corona.  A sequence $(x_1,x_2, \ldots)$ is then a consistent
set of instructions for tiling larger and larger regions of the plane,
which is tantamount to the tiling itself.

Consider two tiles $t_1, t_2$ in (possibly different) tilings in $X_T$ to
be equivalent if a patch of the first tiling, containing $t_1$ and its first
$n$ coronas, is identical, up to translation, to a patch of the second 
tiling containing $t_2$ and its first $n$ coronas. By finite local complexity,
there are only finitely many equivalence clases, which we call $n$-collared
tiles. 

For each $n$-collared tile $t_i$, we consider how such a tile can be
placed at the origin.  This merely means picking a point in $t_i$ to
sit at the origin.  Hence the set of all possible instructions for
placing $t_i$ over the origin is just $t_i$ itself!  A patch of a
tiling in which the origin is on the boundary of two or more tiles is
described by points on the boundary of two or more cells, and
these points must be identified.  The branched manifold $K_n$ is then
the union of all the $n$-collared tiles $t_i$, modulo this
identification.  Since we are using $n$-collared tiles, each of the
points being identified carries complete information about the
placement of the tiles touching the origin, plus their first $n-1$
coronas.  This completes the construction. \qed
  
Note that Theorem \ref{4.1} does not require the tiling to be
generated from a substitution. A substitution tiling space (such as
the space of Penrose tilings) may be constructed as an inverse limit
as above, or by Anderson and Putnam's original contruction \cite{AP}.
Although the construction of Theorem \ref{4.1} is conceptually simple,
it does not lend itself to explicit computations, and the cohomology
of substitution tiling spaces is almost always done using the
Anderson-Putnam construction.

\section{Deformations of tilings}

Since we assumed our tiles to be polygons, describing
the shape and size of a tile means specifying the displacement vector
corresponding to each edge of the polygon.  For instance, the shape of
a triangle is determined by three vectors that add up to zero. 
If, somewhere in a tiling, tiles
$t_1$ and $t_2$ meet along a common edge, then the edge of $t_1$ and 
the matching edge of $t_2$ must have exactly the same displacement vector.  
Furthermore, the sum of all displacement vectors for the edges of $t_i$ must
add up to zero.

That is, the shapes of the tiles are determined by a function $f$ from
the set of tile edges, modulo certain identifications, to $\R^2$. But
the set of tiles, modulo these identifications, is precisely the 
approximant $K_0$!  In other words, a shape is a vector-valued 1-cochain
on $K_0$.  The fact that the sum of the vectors around a tile add to zero
means that this is a {\em closed} cochain.  

We can consider slightly more general deformations, where the new shape of a
tile depends not only on its type but on the types of its nearest neighbours.
In that case the shape is a vector-valued closed 1-cochain on $K_1$.  For 
maximum generality, we may allow closed 1-cochains on any approximant
$K_n$. 

Not all deformations are interesting.  If the 1-cochain is $d$ of a
0-cochain on $K_n$, then the resulting tiling space is MLD to the original
one \cite{CS}.  Our deformations, up to MLD equivalence, are then described by 
closed 1-cochains modulo exact 1-cochains, i.e., by the first cohomology of
$K_n$ with values in $\R^2$.  Since we can take $n$ arbitrarily large, this
means we want an element of the direct limit.  But that's precisely the
cohomology of $X_T$ itself. In other words,

\begin{thm}[CS]
Deformations to a tiling $T$ are parametrized by closed vector-valued
1-cochains on approximants to $X_T$.  If two shape functions define the
same class in $H^1(X_T, R^2)$, then the two tilings are MLD.
\end{thm}

Topological conjugacy is subtler.  Some cohomology classes are
{\em asymptotically negligible}, meaning they represent deformations that 
do not change the topological
conjugacy class of a tiling space. For general tiling spaces these classes
can be difficult to compute, but for substitution tilings (such as Penrose)
they are easy.  The substitution defines a map $\psi: X_T \to X_T$. 
This induces a pullback map $\psi^* H^1(X_T, \R^2) \to H^1(X_T, \R^2)$. 
The eigenspaces of $\psi^*$ whose eigenvalues have magnitude 
strictly less than 1 are asymptotically negligible. [CS]

For the Penrose tiling, $H^1(X_T,\R^2)$ is 10-dimensional.  The eigenvalues
of $\psi^*$ are the golden mean $\tau = (1+ \sqrt{5})/2$, with multiplicity
4, $1-\tau$ with multiplicity 4, and $-1$ with multiplicity 2.  The 
4-dimensional $\tau$ eigenspace describe the 4-dimensional space of linear
transformations (rotations, stretches and shears) that could be applied
uniformly to the tiling.  The $-1$ eigenspace breaks the (statistical)
180-degree rotational symmetry of the Penrose system, and the $1-\tau$
eigenspace is asymptotically negligible.  Since the second tiling in 
Figure 1 preserves the 180-degree rotational symmetry, its shape must be
described by a combination of  $\tau$ and $1-\tau$ eigenenvectors, while
the shape of the original Penrose tiling is entirely a $\tau$ eigenvector. 
Thus the second tiling space is topologically conjugate (but not MLD!)
to a linear transformation applied to the first. 
 
Topological  conjugacies preserve the dynamical spectrum, and therefore
preserve locations (but not necessarily intensities) of the peaks
of the diffraction pattern. As a result, the locations of the
Bragg peaks of the second tiling are related to the Bragg peaks of the first
tiling by a linear transformation.  (For an investigation on how shape
changes affect the intensities of the peaks of the 
1-dimensional Fibonacci tiling, see \cite{Polish}).

\noindent{\em Acknowledgements.} This work was partially supported by
the National Science Foundation.

\end{document}